\documentclass{mmn_art}

\theoremstyle{plain}
\newtheorem{theorem}{Theorem}
\newtheorem{corollary}{Corollary}
\newtheorem{lemma}{Lemma}

\theoremstyle{definition}
\newtheorem{definition}{Definition}

\theoremstyle{remark}
\newtheorem{remark}{Remark}

\numberwithin{equation}{section}

\allowdisplaybreaks[4]

\begin{document}
% Do NOT capitalize the title, section headings and author names, this will be done by our class file!
% Please, provide an abbreviated title.
\title[Inverse problem for a hyperbolic equation with nonlocal conditions]{Nonlinear inverse problem for identifying \\ a coefficient of the lowest term in hyperbolic equation with nonlocal conditions} 

% If there are more than one authors of the manuscript, one author has to be designated as the corresponding one, using the optional parameter "\corresponding" in the command "\address". 

\author{G.~Yu.~Mehdiyeva}
\address{Baku State University \\ Department of Computational Mathematics \\ 23 Z.Khalilov Str.\\ Az1148\\ Baku \\ Azerbaijan}

\email{imn\_bsu@mail.ru}

\author{Y.~T.~Mehraliyev}
\address{Baku State University \\ Department of Differential and Integral Equations \\ 23 Z.Khalilov Str.\\ Az1148\\ Baku \\ Azerbaijan}

\email{yashar\_aze@mail.ru}

%\thanks{The first author was supported in part by the XXX Fund, Grant No.~YYYYY.}

\author{E.~I.~Azizbayov}

\address[\corresponding]{The Academy of Public Administration under the President of the Republic of Azerbaijan \\ Department of Control for Intelligent Systems \\ 74, Lermontov Str. \\ Az 1001 \\  Baku \\ Azerbaijan}

\curraddr{The Academy of Public Administration under the President of the Republic of Azerbaijan \\ Department of Control for Intelligent Systems \\ 74, Lermontov Str. \\ Az 1001 \\  Baku \\ Azerbaijan}

\email{eazizbayov@bsu.edu.az}

%Please, provide an abstract of maximum 120 words. The abstract should state briefly the purpose of the research and its principal results. The usage of references and long math expressions should be avoided.

\begin{abstract}
In this paper, a nonlinear inverse boundary value problem for the second-order hyperbolic equation with nonlocal conditions is studied. To investigate the solvability of the original problem, we first consider an auxiliary inverse boundary value problem and prove its equivalence (in a certain sense) to the original problem. Then using the Fourier method and contraction mappings principle, the existence and uniqueness theorem for auxiliary problem is proved. Further, on the basis of the equivalency of these problems the existence and uniqueness theorem for the classical solution of the considered inverse coefficient problem is proved for the smaller value of time.
\end{abstract}

%\dedicatory{This paper is dedicated to Professor X on his 125th birthday.}

%Please, provide the 2010 Mathematics Subject Classification code(s)

\subjclass{35R30; 35L10; 35A01; 35A02; 35A09}

% Please, provide a maximum of 6 keywords, separated by comma. Do NOT capitalize keywords!
% Keywords are short phrases, which will be used for indexing purposes.

\keywords{inverse problem, hyperbolic equation, overdetermination condition, classical solution, existence, uniqueness}

\pdfmetadata

\maketitle

\section{Introduction and problem statement}

Let $0 < T <  + \infty$ be some fixed number and $D_T$ be a rectangular region in the $xt$-plane defined by $0 \le x \le 1, \ 0 \le t \le T$. Consider the problem of determining the unknown functions $u = u(x,t)$ and $a = a(t)$ such that the pair $\{u,a\}$ satisfies the following hyperbolic equation of second order
\begin{equation}\label{E:1}
u_{tt} (x,t) - u_{xx} (x,t) = a(t)u(x,t) + f(x,t) \ \ (x,t) \in D_T ,
\end{equation}
with the nonlocal initial conditions
\begin{equation}\label{E:2}
u(x,0) + \delta _1 u(x,T) = \varphi (x), \ \ u_t (x,0) + \delta _2 u_t (x,T) = \psi (x), \ 0 \le x \le 1,
\end{equation}
the boundary conditions
\begin{equation}\label{E:3}
u(0,t) = \beta u(1,t), \ 0 \le t \le T,
\end{equation}
\begin{equation}\label{E:4}
\int\limits_0^1 {u(x,t)dx}  = 0, \ 0 \le t \le T,
\end{equation}
and overdetermination condition
\begin{equation}\label{E:5}
u\left( {\frac{1}{2},t} \right) = h(t), \ 0 \le t \le T,
\end{equation}
in which $\delta _1 , \delta _2  \ge 0$ and $\beta  \ne  \pm 1$ are given numbers, $f(x,t),\varphi (x)$, and $\psi (x)$ are known functions of $x \in [0,1]$ and $t \in [0,T]$.

In the present work we investigate a nonlinear inverse problem for identifying a coefficient of the lowest term in hyperbolic equation from the overdetermination data. Such problems are called inverse problems in mathematical physics. The applied importance of inverse problems is so great (seismology, mineral exploration, biology, medicine, desalination of seawater, movement of liquid in a porous medium, acoustics, electromagnetics, fluid dynamics, calculating the density of the Earth from measurements of its gravity field, for example) which puts them a series of the most actual problems of modern mathematics. Strictly speaking, the inverse problems for hyperbolic/wave equations are of prime interest in seismology. Besides, vibrations of structures (as buildings and beams) are modeled by hyperbolic differential equations. However, it should be noted that an inverse problem is called linear if the recovery function enters the given equation linearly, and nonlinear otherwise.

Nowadays, inverse problems for hyperbolic equations have been well studied by many authors using different methods (see, e.g., \cite{AlLesDin,CanDun,Eskin,JiangYam,KozSaf,BaJiLi,Salazar,Shcheglov,SloSel,UhlZha}, and the references given therein). Moreover, in \cite{Rag01,Rag12,RagScape17,RagTach20} the authors present a regularity result for solutions of partial differential equations in the framework of mixed Morrey spaces.

A distinctive feature of presented article is the investigation of an inverse hyperbolic problem with both spatial and time nonlocal conditions.

This article is based on ideas close to those used in \cite{AlMeh,Azizbayov,MehIsg,TekMehIs}.

\begin{definition}
A pair of functions $\{ u(x,t),a(t)\}$ is said to be a classical solution of problem \eqref{E:1}--\eqref{E:5} if all three of the following conditions are satisfied:
\begin{itemize}
 \item [{a.}]
The function $u(x,t)$ with the derivatives $u_{xx} (x,t)$ and $u_{tt} (x,t)$ are continuous in the domain $D_T$.
 \item [{b.}]
The function $a(t)$ is continuous on the interval $[0,T]$.
\item [{c.}]
The Eq. \eqref{E:1} and conditions \eqref{E:2}--\eqref{E:5} are satisfied in the classical (usual) sense.
\end{itemize}
\end{definition}

Now, to study problem \eqref{E:1}--\eqref{E:5}, we consider the following auxiliary inverse boundary value problem: it is required to find a pair of functions $u(x,t) \in C^2 (D_T ),\\ \ a(t) \in C[0,T]$ from \eqref{E:1}--\eqref{E:3} and
\begin{equation}\label{E:6}
u_x (0,t) = u_x (1,t), \ 0 \le t \le T,
\end{equation}
\begin{equation}\label{E:7}
h''(t) - u_{xx} \left( {\frac{1}{2},t} \right) = a(t)h(t) + f\left( {\frac{1}{2},t} \right), \ 0 \le t \le T.
\end{equation}

Similarly (see \cite{MehAz}, Theorem 2.2, p.4) it can be proved that

\begin{theorem} \label{thm:1}
Suppose that $\varphi (x), \psi (x) \in C[0,1], \ h(t) \in C^2 [0,T], \ h(t) \ne 0$, \newline $ f(x,t) \in C(D_T ), \ \int\limits_0^1 {f(x,t)dx}  = 0, \ 0 \le t \le T$ and the compatibility conditions
\[
\int\limits_0^1 {\varphi (x)dx}  = 0, \ \ \int\limits_0^1 {\psi (x)dx}  = 0,
\]
\[
h(0) + \delta _1 h(T) = \varphi \left( {\frac{1}{2}} \right), \ \ h'(0) + \delta _2 h'(T) = \psi \left( {\frac{1}{2}} \right),
\]
hold. Then the following statements are true:
\begin{itemize}
 \item [(\emph{i})]
Each classical solution $\{u(x,t),a(t)\}$ of problem \eqref{E:1}--\eqref{E:5} is the solution of problem \eqref{E:1}--\eqref{E:3}, \eqref{E:6}, \eqref{E:7}, as well
\item [(\emph{ii})]
each solution $\{u(x,t),a(t)\}$ of problem \eqref{E:1}--\eqref{E:3}, \eqref{E:6}, \eqref{E:7} is a classical solution to the problem \eqref{E:1}--\eqref{E:5}, if
\[
\frac{{(1 + 2\delta _1  + 3\delta _2  + \delta _1 \delta _2 )T^2 }}{{2(1 + \delta _1 )(1 + \delta _2 )}}\left\| {a(t)} \right\|_{C[0,T]}  < 1.
\]
\end{itemize}
\end{theorem}

\section{Auxiliary facts and denotations}\label{Sec:Auxiliary facts and denotations}
It is known that sequences of functions \cite{MehYus}
\begin{equation}\label{E:8}
 X_0 (x) = px + q, \ \ X_{2k - 1} (x) = (px + q)\cos \lambda _k x, \ \ X_{2k} (x) = \sin \lambda _k x, \ \ k = 1,2, \ldots,
\end{equation}
\begin{equation}\label{E:9}
Y_0 (x) = 2, \ \ Y_{2k - 1} (x) = 4\sin \lambda _k x, \ \ Y_{2k} (x) = q(1 - q - px)\cos \lambda _k x, \ \ k = 1,2, \ldots,
\end{equation}
form a biorthogonal system and system \eqref{E:8} forms a Riesz basis in $L_2 (0,1)$ for $\lambda_k  = 2k\pi$ \ $(k = 1,2, \ldots)$. Here $p$ and $q$ denotes, in turn, the numbers
\[
p = \frac{{1 - \beta }}{{1 + \beta }} \ne 0, \ \ q = \frac{\beta }{{1 + \beta }}.
\]

We state the following lemmas without proof.

\begin{lemma}\label{lem:1}
(see {\cite{AzMeh,Mehraliyev}})\label{L:2.1} For any function $v(x)$ with the properties:
\[
v(x) \in C^{2i - 1} [0,1], \ \ v^{(2i)} (x) \in L_2 (0,1),
\]
\[
 v^{(2s)} (0) = \beta v^{(2s)} (1), \ \  v^{(2s + 1)} (0) = v^{(2s + 1)} (1) \ \ (i \ge 1; \ \ s = \overline {0,i}),
\]
the estimates are valid
\[
\left( {\sum\limits_{k = 0}^\infty  {(\lambda _k^{2i} v_{2k - 1} )^2 } } \right)^{\frac{1}{2}}  \le 2\sqrt 2 \left\| {v^{(2i)} (x)} \right\|_{L_2 (0,1)} ,
\]
\begin{equation}\label{E:10}
\left( {\sum\limits_{k = 1}^\infty  {(\lambda _k^{2i} v_{2k} )^2 } } \right)^{\frac{1}{2}}\le 2\sqrt 2 \left\| {v^{(2i)} (x)(1 - q - px) - 2ipv^{(2i - 1)} (x)} \right\|_{L_2 (0,1)} ,
\end{equation}
where
\[
v_k  = \int\limits_0^1 {v(x)Y_k (x)dx}, \ k = 0,1,....
\]
\end{lemma}

\begin{lemma}\label{lem:2}
(see {\cite{AzMeh,Mehraliyev}}) Under the assumptions:
\[
v(x) \in C^{2i} [0,1], \ \ v^{(2i + 1)} (x) \in L_2 (0,1),
\]
\[ v^{(2s)} (0) = \beta v^{(2s)} (1), \ \ v^{(2s + 1)} (0) = v^{(2s + 1)} (1) \ \ (i \ge 1; \ \ s = \overline {0,i}),
\]
we establish the validity of the estimates
\[
\left( {\sum\limits_{k = 1}^\infty  {(\lambda _k^{2i + 1} v_{2k - 1} )^2 } } \right)^{\frac{1}{2}}\le 2\sqrt 2 \left\| {v^{(2i + 1)} (x)} \right\|_{L_2 (0,1)} ,
\]
\[
\left( {\sum\limits_{k = 1}^\infty  {(\lambda _k^{2i + 1} v_{2k - 1} )^2 } } \right)^{\frac{1}{2}}
\]
\begin{equation}\label{E:11}
\le 2\sqrt 2 \left\| {v^{(2i + 1)} (x)(1 - q - px) - (2i + 1)pv^{(2i)} (x)} \right\|_{L_2 (0,1)} .
\end{equation}
\end{lemma}

We now look at the following functional spaces:

$B_{2,T}^3$ \cite{Mehraliyev} denotes a set of all functions of the form
\[
u(x,t) = \sum\limits_{k = 1}^\infty  {u_k (t)X_k (x)} ,
\]
considered in $D_T$, where $u_k (t) \in C[0,T]$ and
\[
J(u) \equiv \left\| {u_0 (t)} \right\|_{C[0,T]}
\]
\[
+ \left\{ {\sum\limits_{k = 1}^\infty  {(\lambda _k^3 \left\| {u_{2k - 1} (t)} \right\|_{C[0,T]} )^2 } } \right\}^{\frac{1}{2}}  + \left\{ {\sum\limits_{k = 1}^\infty  {(\lambda _k^3 \left\| {u_{2k} (t)} \right\|_{C[0,T]} )^2 } } \right\}^{\frac{1}{2}}  <  + \infty .
\]

Actually, the functions $X_k (x) \ (k = 0,1, \ldots )$ defined by the relation \eqref{E:8}. The norm on the set $J(u)$ is established as follows
\[
\left\| {u(x,t)} \right\|_{B_{2,T}^3 }  = J(u).
\]

Let $E_T^3$ denote the space consisting of the topological product $B_{2,T}^3 \times C[0,T]$, which is the norm of the element $z = \{ u,a\}$ defined by the formula
\[
\left\| {z(x,t)} \right\|_{E_T^3 }  = \left\| {u(x,t)} \right\|_{B_{2,T}^3 }  + \left\| {a(t)} \right\|_{C[0,T]} .
\]

\section{Classical solvability of inverse boundary value problem}\label{Sec:Classical solvability of inverse boundary value problem}

Suppose that the data of problem \eqref{E:1}--\eqref{E:3},\eqref{E:6},\eqref{E:7} satisfies the following conditions:

\begin{itemize}
\item [$C_1$.]
$
\delta _1 ,\delta _2  \ge 0, \ \ 1 + \delta _1 \delta _2  \ge \delta _1  + \delta _2 ;
$
\item [$C_2$.]
$
\varphi (x) \in C^2 [0,1], \ \varphi '''(x) \in L_2 (0,1), \ \varphi (0) = \beta \varphi (1), \\ \varphi '(0) = \varphi '(1), \ \varphi ''(0) = \beta \varphi ''(1);
$
\item [$C_3$.]
$
\psi (x) \in C^1 [0,1], \ \psi ''(x) \in L_2 (0,1), \ \psi (0) = \beta \psi (1), \ \psi '(0) = \psi '(1);
$
\item [$C_4$.]
$
f(x,t),f_x (x,t) \in C(D_T ),\ f_{xx} (x,t) \in L_2 (D_T ), \ f(0,t) = \beta f(1,t), \\ f_x (0,t) = f_x (1,t), \ 0 \le t \le T;
$
\item [$C_5$.]
$
h(t) \in C^2 [0,T], \ h(t) \ne 0, \ 0 \le t \le T.
$
\end{itemize}

Since the system \eqref{E:8} forms a Riesz basis in $L_2 (0,1)$. Then each solution to problem \eqref{E:1}--\eqref{E:3},\eqref{E:6},\eqref{E:7} can be sought in the form:
\begin{equation}\label{E:12}
u(x,t) = \sum\limits_{k = 0}^\infty  {u_k (t)X_k (x)} ,
\end{equation}
where
\[
u_k (t) = \int\limits_0^1 {u(x,t)Y_k (x)dx} , \ k = 0,1,....
\]

Moreover, $X_k (x) \ (k = 0,1,...)$ and $Y_k (x) \ (k = 0,1,...)$ are defined by relations \eqref{E:8} and \eqref{E:9}, respectively.

To determine of the desired functions $u_k (t) \ (k = 0,1,2,...)$, using separation of variables and by \eqref{E:1} and \eqref{E:2}, we get
\begin{equation}\label{E:13}
u''_0 (t) = F_0 (t;u,a), \ 0 \le t \le T,
\end{equation}
\begin{equation}\label{E:14}
u''_{2k - 1} (t) + \lambda _k^2 u_{2k - 1} (t) = F_{2k - 1} (t;u,a), \ k = 1,2,...; \ 0 \le t \le T,
\end{equation}
\begin{equation}\label{E:15}
u''_{2k} (t) + \lambda _k^2 u_{2k} (t) = F_{2k} (t;u,a) - 2p\lambda _k u_{2k - 1} (t), \ k = 1,2,...; \ 0 \le t \le T,
\end{equation}
\begin{equation}\label{E:16}
u_k (0) + \delta _1 u_k (T) = \varphi _k , \ \ u'_k (0) + \delta _2 u'_k (T) = \psi _k , \ k = 0,1,2, \ldots ,
\end{equation}
where
\[
F_k (t;u,a) = a(t)u_k (t) + f_k (t), \ \ f_k (t) = \int\limits_0^1 {f(x,t)Y_k (x)dx} , \ k = 0,1,....
\]
\[
\varphi _k  = \int\limits_0^1 {\varphi (x)Y_k (x)dx} , \ \ \psi _k  = \int\limits_0^1 {\psi (x)Y_k (x)dx} , \ k = 0,1,2....
\]

Solving problem \eqref{E:13}--\eqref{E:16}, we find
\begin{equation}\label{E:17}
u_0 (t) = \frac{{\varphi _0 }}{{1 + \delta _1 }} + \frac{{t - \delta _1 (T - t)}}{{(1 + \delta _1 )(1 + \delta _2 )}}\psi _0  + \int\limits_0^T {G_0 (t,\tau )F_0 (\tau ;u,a)d\tau } ,
\end{equation}
\[
u_{2k - 1} (t) = \frac{1}{{\rho _k (T)}}\left[ {\varphi _{2k - 1} (\cos \lambda _k t + \delta _2 \cos \lambda _k (T - t)) } \right.
\]
\begin{equation}\label{E:18}
+ \left. {\frac{{\psi _{2k - 1} }}{{\lambda _k }}(\sin \lambda _k t - \delta _1 \sin \lambda _k (T - t))} \right] + \int\limits_0^T {G_k (t,\tau )F_{2k - 1} (\tau ;u,a)d\tau } ,
\end{equation}
\[
u_{2k} (t) = \frac{1}{{\rho _k (T)}}\left[ {\varphi _{2k} (\cos \lambda _k t + \delta _2 \cos \lambda _k (T - t)) + \frac{{\psi _{2k} }}{{\lambda _k }}(\sin \lambda _k t - \delta _1 \sin \lambda _k (T - t))} \right]
\]
\[
 + \int\limits_0^T {G_k (t,\tau )(F_{2k} (\tau ;u,a))d\tau }
  \]
  \[
-\varphi _{2k - 1} \left\{ { - \frac{1}{{\rho _k^2 (T)}}\frac{1}{{\lambda _k }}} \right.\left[ {\delta _1 \cos \lambda _k t\left( {\frac{T}{2}\sin \lambda _k T + \delta _2 \frac{1}{4}(1 - \cos 2\lambda _k T)} \right)} \right.
\]
\[
 + \delta _2 \sin \lambda _k t\left( {\frac{1}{{2\lambda _k }}\sin \lambda _k T + \frac{T}{2}\cos \lambda _k T + \delta _2 \left( {\frac{T}{2} + \frac{1}{{4\lambda _k }}\sin 2\lambda _k T} \right)} \right)
\]
\[
 + \delta _1 \delta _2 \left( {\frac{1}{{4\lambda _k }}(\cos \lambda _k (2T - t) - \cos \lambda _k t) + \frac{T}{2}\sin \lambda _k t} \right.
\]
\[
\left. {\left. { + \delta _2 \left( { - \frac{T}{2}\sin \lambda _k (T - t) + \frac{1}{{4\lambda _k }}(\cos \lambda _k (T - t) - \cos \lambda _k (T + t))} \right)} \right)} \right]
\]
\[
\left. { + \frac{t}{2}\sin \lambda _k t + \delta _2 \left( { - \frac{t}{2}\sin \lambda _k (T - t) + \frac{1}{{4\lambda _k }}(\cos \lambda _k (T - t) - \cos \lambda _k (T + t))} \right)} \right\}
\]
\[
 + \frac{{\psi _{2k - 1} }}{{\lambda _k }}\left\{ { - \frac{1}{{\rho _k^2 (T)}} \cdot \frac{1}{{\lambda _k }}} \right.\left[ {\delta _1 \cos \lambda _k t\left( {\frac{1}{{2\lambda _k }}\sin \lambda _k T - \frac{T}{2}\cos \lambda _k T} \right)} \right.
\]
\[
 - \delta _1 \left( {\frac{T}{2} - \frac{1}{{4\lambda _k }}\sin 2\lambda _k T} \right) + \delta _2 \sin \lambda _k t\left( {\frac{T}{2}\sin \lambda _k T - \frac{{\delta _1 }}{{4\lambda _k }}(1 - \cos 2\lambda _k T)} \right)
\]
\[
 + \delta _1 \delta _2 \left( {\frac{1}{{4\lambda _k }}(\sin \lambda _k (2T - t) + \sin \lambda _k t) - \frac{T}{2}\cos \lambda _k t} \right.
\]
\[
\left. {\left. { - \delta _1 \left( {\frac{T}{2}\cos \lambda _k (T - t) - \frac{1}{{4\lambda _k }}(\sin \lambda _k (T - t) + \sin \lambda _k (T + t))} \right)} \right)} \right]
\]
\[
 + \left( {\frac{1}{{2\lambda _k }}\sin \lambda _k t - \frac{t}{2}\cos \lambda _k t - \delta _1 \left( {\frac{t}{2}\cos \lambda _k (T - t) + } \right.} \right.\frac{1}{{4\lambda _k }}(\sin \lambda _k (T - t)
\]
\begin{equation}\label{E:19}
\left. {\left. { - \sin \lambda _k (T + t))} \right)} \right\} + \int\limits_0^T {G_k (t,\tau )\left( {\int\limits_0^T {G_k (\tau ,\xi )F_{2k - 1} (\xi ;u,a)d\xi } } \right)} d\tau ,
\end{equation}
where
\begin{equation*}
	G_0 (t,\tau ) =
	\begin{cases}
		 - \frac{{\delta _2 t + \delta _1 (T - \tau ) + \delta _1 \delta _2 (t - \tau )}}{{(1 + \delta _1 )(1 + \delta _2 )}},\ \ \  \ \ \ \ \ \ \ \ t \in [0,\tau ], \\
		- \frac{{\delta _2 t + \delta _1 (T - \tau ) - (1 + \delta _1  + \delta _2 )(t - \tau )}}{{(1 + \delta _1 )(1 + \delta _2 )}}, \ \ \ \ t \in [\tau ,T], 
	\end{cases}
\end{equation*}

\begin{equation*}
	G_k (t,\tau ) =
	\begin{cases}
		 - \frac{1}{{\rho _k (T)}} \cdot \frac{1}{{\lambda _k }}[\delta _1 \sin \lambda _k (T - \tau )\cos \lambda _k t + \delta _2 \cos \lambda _k (T - \tau )\sin \lambda _k t
 \\ \ \ \ \ \ \ \ \ \ \ \ \ \ \ \ \ \ \ \  \ \ \ \ \ \ \ \ \ \ \ \ \ \ \ \ \ \ \ \ \ \ \ \ \  \ \ \ \ \ \ \ \  \ \ \  
 + \delta _1 \delta _2 \sin \lambda _k (t - \tau )], \ t \in [0,\tau ], \\
 - \frac{1}{{\rho _k (T)}} \cdot \frac{1}{{\lambda _k }}[\delta _1 \sin \lambda _k (T - \tau )\cos \lambda _k t + \delta _2 \cos \lambda _k (T - \tau )\sin \lambda _k t
 \\ \ \ \ \ \ \ \ \ \ \ \ \ \ \ \ \ \ \ \  \ \ \ \  + \delta _1 \delta _2 \sin \lambda _k (t - \tau )] + \frac{1}{{\lambda _k }}\sin \lambda _k (t - \tau ),\ t \in [\tau ,T], \\
	\end{cases}
\end{equation*}

and
\[
\rho _k (T) = 1 + (\delta _1  + \delta _2 )\cos \lambda _k T + \delta _1 \delta _2 .
\]

Substituting the expressions of \eqref{E:17}, \eqref{E:18}, and \eqref{E:19} into \eqref{E:12}, we find the component $u(x,t)$ of the classical solution to problem \eqref{E:1}--\eqref{E:3}, \eqref{E:6}, \eqref{E:7} to be
\[
u(x,t) = \left\{ {\frac{{\varphi _0 }}{{1 + \delta _1 }} + \frac{{t - \delta _1 (T - t)}}{{(1 + \delta _1 )(1 + \delta _2 )}}\psi _0  + \int\limits_0^T {G_0 (t,\tau )F_0 (\tau ;u,a)d\tau } } \right\}X_0 (x)
\]
\[
+ \sum\limits_{k = 1}^\infty  {\left\{ {\frac{1}{{\rho _k (T)}}\left[ {\varphi _{2k - 1} (\cos \lambda _k t + \delta _2 \cos \lambda _k (T - t))} \right.} \right.} 
\]
\[
\left. {\left. { + \frac{{\psi _{2k - 1} }}{{\lambda _k }}(\sin \lambda _k t - \delta _1 \sin \lambda _k (T - t))} \right] + \int\limits_0^T {G_k (t,\tau )F_{2k - 1} (\tau ;u,a)d\tau } } \right\}X_{2k - 1} (x)
\]
\[
 + \sum\limits_{k = 1}^\infty  {\left\{ {\frac{1}{{\rho _k (T)}}[\varphi _{2k} (\cos \lambda _k t + \delta _2 \cos \lambda _k (T - t))+ } \right.} 
\]
\[
 + \left. {\frac{{\psi _{2k} }}{{\lambda _k }}(\sin \lambda _k t - \delta _1 \sin \lambda _k (T - t))} \right] + \int\limits_0^T {G_k (t,\tau )F_{2k} (\tau ;u,a)d\tau }
\]
\[
 - \varphi _{2k - 1} \left\{ { - \frac{1}{{\rho _k^2 (T)}} \cdot \frac{1}{{\lambda _k }}} \right.\left[ {\delta _1 \cos \lambda _k t\left( {\frac{T}{2}\sin \lambda _k T + \delta _2 \frac{1}{4}(1 - \cos 2\lambda _k T)} \right)} \right.
\]
\[
 + \delta _2 \sin \lambda _k t\left( {\frac{1}{{2\lambda _k }}\sin \lambda _k T + \frac{T}{2}\cos \lambda _k T + \delta _2 \left( {\frac{T}{2} + \frac{1}{{4\lambda _k }}\sin 2\lambda _k T} \right)} \right)
\]
\[
 + \delta _1 \delta _2 \left( {\frac{1}{{4\lambda _k }}(\cos \lambda _k (2T - t) - \cos \lambda _k t) + \frac{T}{2}\sin \lambda _k t} \right.
\]
\[
\left. {\left. { + \delta _2 \left( { - \frac{T}{2}\sin \lambda _k (T - t) + \frac{1}{{4\lambda _k }}(\cos \lambda _k (T - t) - \cos \lambda _k (T + t))} \right)} \right)} \right]
\]
\[
\left. { + \frac{t}{2}\sin \lambda _k t + \delta _2 \left( { - \frac{t}{2}\sin \lambda _k (T - t) + \frac{1}{{4\lambda _k }}(\cos \lambda _k (T - t) - \cos \lambda _k (T + t))} \right)} \right\}
\]
\[
 + \frac{{\psi _{2k - 1} }}{{\lambda _k }}\left\{ { - \frac{1}{{\rho _k^2 (T)}} \cdot \frac{1}{{\lambda _k }}} \right.\left[ {\delta _1 \cos \lambda _k t\left( {\frac{1}{{2\lambda _k }}\sin \lambda _k T - \frac{T}{2}\cos \lambda _k T} \right)} \right.
\]
\[
 - \delta _1 \left( {\frac{T}{2} - \frac{1}{{4\lambda _k }}\sin 2\lambda _k T} \right) + \delta _2 \sin \lambda _k t\left( {\frac{T}{2}\sin \lambda _k T - \frac{{\delta _1 }}{{4\lambda _k }}(1 - \cos 2\lambda _k T)} \right)
\]
\[
 + \delta _1 \delta _2 \left( {\frac{1}{{4\lambda _k }}(\sin \lambda _k (2T - t) + \sin \lambda _k t) - \frac{T}{2}\cos \lambda _k t} \right.
\]
\[
\left. {\left. { - \delta _1 \left( {\frac{T}{2}\cos \lambda _k (T - t) - \frac{1}{{4\lambda _k }}(\sin \lambda _k (T - t) + \sin \lambda _k (T + t))} \right)} \right)} \right]
\]
\[
+ \left( {\frac{1}{{2\lambda _k }}\sin \lambda _k t - \frac{t}{2}\cos \lambda _k t - \delta _1 \left( {\frac{t}{2}\cos \lambda _k (T - t) + } \right.} \right.\frac{1}{{4\lambda _k }}(\sin \lambda _k (T - t)
\]
\begin{equation}\label{E:20}
- \sin \lambda _k (T + t)))\}  + \left. {\int\limits_0^T {G_k (t,\tau )\left( {\int\limits_0^T {G_k (\tau ,\xi )F_{2k - 1} (\xi ;u,a)d\xi } } \right)} d\tau } \right\}X_{2k} (x).
\end{equation}

Now, using \eqref{E:12} from \eqref{E:7} we have
\begin{equation}\label{E:21}
a(t) = [h(t)]^{ - 1} \left\{ {h''(t) - f\left( {\frac{1}{2},t} \right)} \right. - \frac{1}{2}\left. {\sum\limits_{k = 1}^\infty  {( - 1)^k \lambda _k^2 } u_{2k - 1} (t)} \right\}.
\end{equation}

Substituting expressions $u_{2k - 1} (t) \ (k = 1,2,...)$ from \eqref{E:19} into \eqref{E:21}, immediately yields:
\[
a(t) = [h(t)]^{ - 1}
\]
\[\times \left\{ {h''(t) - f\left( {\frac{1}{2},t} \right)} \right. - \frac{1}{2}\sum\limits_{k = 1}^\infty  {( - 1)^k \lambda _k^2 \left\{ {\frac{1}{{\rho _k (T)}}\left[ {\varphi _{2k - 1} (\cos \lambda _k t + \delta _2 \cos \lambda _k (T - t)) } \right.} \right.}
\]
\begin{equation}\label{E:22}
+ \left. {\left. {\left. {\frac{{\psi _{2k - 1} }}{{\lambda _k }}(\sin \lambda _k t - \delta _1 \sin \lambda _k (T - t))} \right] + \int\limits_0^T {G_k (t,\tau )F_{2k - 1} (\tau ;u,a)d\tau } } \right\}} \right\}.
\end{equation}

Thus, the solution of problem \eqref{E:1}--\eqref{E:3}, \eqref{E:6}, \eqref{E:7} is reduced to the solution of system \eqref{E:20}, \eqref{E:22} with respect to unknown functions $u(x,t)$ and $a(t)$.

Using the same discussion in \cite{Mehraliyev}, one can prove

\begin{lemma}\label{lem:3}
If $\{ u(x,t),a(t)\}$ is any solution to problem \eqref{E:1}--\eqref{E:3}, \eqref{E:6}, \eqref{E:7}, then the functions
\[
u_k (t) = \int\limits_0^1 {u(x,t)Y_k (x)dx} , \ \ k = 0,1,...,
\]
satisfy the countable system \eqref{E:17}, \eqref{E:18} and \eqref{E:19} on an interval $[0,T]$.
\end{lemma}

Obviously, if $u_k (t) = \int\limits_0^1 {u(x,t)Y_k (x)dx} , \ k = 0,1,...$ is a solution to system \eqref{E:17}, \eqref{E:18} and \eqref{E:19}, then the pair $\{ u(x,t),a(t)\}$ of functions $u(x,t) = \sum\limits_{k = 0}^\infty  {u_k (t)} X_k (x)$ and $a(t)$ is also a solution to system \eqref{E:20}, \eqref{E:22}.

It follows from the Lemma 3 that

\begin{corollary}
 Let system \eqref{E:20}, \eqref{E:22} have a unique solution. Then problem \eqref{E:1}--\eqref{E:3}, \eqref{E:6}, \eqref{E:7} cannot have more than one solution, i.e. if the problem \eqref{E:1}--\eqref{E:3}, \eqref{E:6}, \eqref{E:7} has a solution, then it is unique.
\end{corollary}

Now consider the operator
\[
\Phi (u,a) = \{ \Phi _1 (u,a),\Phi _2 (u,a)\} ,
\]
in the space $E_T^3$, where
\[
\Phi _1 (u,a) = \tilde u(x,t) \equiv \tilde u_0 (t)X_0 (x) + \sum\limits_{k = 0}^\infty  {\tilde u_{2k - 1} (t)X_{2k - 1} (x)}  + \sum\limits_{k = 1}^\infty  {\tilde u_{2k} (t)X_{2k} (x)},
\]
\[
\Phi _2 (u,a) = \tilde a(t),
\]
and the functions $\tilde u_0 (t),\tilde u_{2k - 1} (t),\tilde u_{2k} (t) \ (k = 1,2,...)$, and $\tilde a(t)$ are equal correspondingly to the right sides of \eqref{E:17}, \eqref{E:18}, \eqref{E:19}, and \eqref{E:22}.

It is clear that under conditions $\delta _1 ,\delta _2  \ge 0, \ \ 1 + \delta _1 \delta _2  > \delta _1  + \delta _2$ and using the inequality
\[
\rho _k (T) \ge 1 - (\delta _1  + \delta _2 ) + \delta _1 \delta _2 
\]
we can write
\[
\frac{1}{{\rho _k (T)}} \le \frac{1}{{1 - (\delta _1  + \delta _2 ) + \delta _1 \delta _2 }} \equiv \rho  > 0.
\]

Using this relation, from \eqref{E:17}, \eqref{E:18}, \eqref{E:19}, and \eqref{E:22} we obtain
\[
\left\| {\tilde u_0 (t)} \right\|_{C[0,T]}  \le \frac{1}{{1 + \delta _1 }}\left| {\varphi _0 } \right| + \frac{T}{{1 + \delta _2 }}\psi _0  + \frac{{1 + 3\delta _1  + 3\delta _2 }}{{(1 + \delta _1 )(1 + \delta _2 )}}T
\]
\begin{equation}\label{E:23}
\times \left( {\sqrt T \left( {\int\limits_0^T {\left| {f_0 (\tau )} \right|^2 d\tau } } \right)^{\frac{1}{2}}  + T\left\| {a(t)} \right\|_{C[0,T]} \left\| {u_0 (t)} \right\|_{C[0,T]} } \right),
\end{equation}
\[
\left\{ {\sum\limits_{k = 1}^\infty  {(\lambda _k^3 \left\| {\tilde u_{2k - 1} (t)} \right\|_{C[0,T]} )^2 } } \right\}^{\frac{1}{2}}
\]
\[
\le 2\rho (1 + \delta _2 )\left( {\sum\limits_{k = 1}^\infty  {(\lambda _k^3 \left| {\varphi _{2k - 1} } \right|)^2 } } \right)^{\frac{1}{2}}  + 2\rho (1 + \delta _1 )\left( {\sum\limits_{k = 1}^\infty  {(\lambda _k^2 \left| {\psi _{2k - 1} } \right|)^2 } } \right)^{\frac{1}{2}}
\]
\[
 + 2(1 + 2\rho (\delta _1  + \delta _2  + \delta _1 \delta _2 ))\sqrt T \left( {\int\limits_0^T {\sum\limits_{k = 1}^\infty  {(\lambda _k^2 \left| {f_{2k - 1} (\tau )} \right|)^2 d\tau } } } \right)^{\frac{1}{2}}
\]
\begin{equation}\label{E:24}
+ 2(1 + 2\rho (\delta _1  + \delta _2  + \delta _1 \delta _2 ))T\left\| {a(t)} \right\|_{C[0,T]} \left( {\sum\limits_{k = 1}^\infty  {(\lambda _k^3 \left\| {u_{2k - 1} (t)} \right\|_{C[0,T]} )^2 } } \right)^{\frac{1}{2}} ,
\end{equation}
\[
\left\{ {\sum\limits_{k = 1}^\infty  {(\lambda _k^3 \left\| {\tilde u_{2k} (t)} \right\|_{C[0,T]} )^2 } } \right\}^{\frac{1}{2}}
\]
\[
\le 2\sqrt 2 \rho (1 + \delta _2 )\left( {\sum\limits_{k = 1}^\infty  {(\lambda _k^3 \left| {\varphi _{2k} } \right|)^2 } } \right)^{\frac{1}{2}}  + 2\sqrt 2 \rho (1 + \delta _1 )\left( {\sum\limits_{k = 1}^\infty  {(\lambda _k^2 \left| {\psi _{2k} } \right|)^2 } } \right)^{\frac{1}{2}}
\]
\[
 + 2\sqrt 2 (1 + 2\rho (\delta _1  + \delta _2  + \delta _1 \delta _2 ))\sqrt T \left( {\int\limits_0^T {\sum\limits_{k = 1}^\infty  {(\lambda _k^2 \left| {f_{2k} (\tau )} \right|)^2 d\tau } } } \right)^{\frac{1}{2}}
\]
\[
+ 2\sqrt 2 (1 + 2\rho (\delta _1  + \delta _2  + \delta _1 \delta _2 ))T\left\| {a(t)} \right\|_{C[0,T]} \left( {\sum\limits_{k = 1}^\infty  {(\lambda _k^3 \left\| {u_{2k} (t)} \right\|_{C[0,T]} )^2 } } \right)^{\frac{1}{2}} \]
\[ + 2\sqrt 2 \rho _1 \left( {\sum\limits_{k = 1}^\infty  {(\lambda _k^3 \left| {\varphi _{2k - 1} } \right|)^2 } } \right)^{\frac{1}{2}} + 2\sqrt 2 \rho _2 \left( {\sum\limits_{k = 1}^\infty  {(\lambda _k^2 \left| {\psi _{2k - 1} } \right|)^2 } } \right)^{\frac{1}{2}}
\]
\[
+ 2\sqrt 2 (1 + 2\rho (\delta _1  + \delta _2  + \delta _1 \delta _2 ))^2 T\left[ {\sqrt T \left( {\int\limits_0^T {\sum\limits_{k = 1}^\infty  {(\lambda _k^2 \left| {f_{2k - 1} (\tau )} \right|)^2 d\tau } } } \right)^{\frac{1}{2}} } \right.
\]
\begin{equation}\label{E:25}
\left. { + T\left\| {a(t)} \right\|_{C[0,T]} \left( {\sum\limits_{k = 1}^\infty  {(\lambda _k^3 \left\| {u_{2k - 1} (t)} \right\|_{C[0,T]} )^2 } } \right)^{\frac{1}{2}} } \right],
\end{equation}
\[
\left\| {\tilde a(t)} \right\|_{C[0,T]}  \le \left\| {[h(t)]^{ - 1} } \right\|_{C[0,T]} \left\{ {\left\| {h''(t) - f\left( {\frac{1}{2},t} \right)} \right\|_{C[0,T]}  + \frac{1}{2}\left( {\sum\limits_{k = 1}^\infty  {\lambda _k^{ - 2} } } \right)^{\frac{1}{2}} } \right.
\]
\[
 \times \left[ {\rho (1 + \delta _2 )\left( {\sum\limits_{k = 1}^\infty  {(\lambda _k^3 \left| {\varphi _{2k - 1} } \right|)^2 } } \right)^{\frac{1}{2}}  + \rho (1 + \delta _1 )\left( {\sum\limits_{k = 1}^\infty  {(\lambda _k^2 \left| {\psi _{2k - 1} } \right|)^2 } } \right)^{\frac{1}{2}} } \right.
\]
\[
 + (1 + 2\rho (\delta _1  + \delta _2  + \delta _1 \delta _2 ))\sqrt T \left( {\int\limits_0^T {\sum\limits_{k = 1}^\infty  {(\lambda _k^2 \left| {f_{2k - 1} (\tau )} \right|)^2 d\tau } } } \right)^{\frac{1}{2}}
\]
\begin{equation}\label{E:26}
\left. {\left. { + (1 + 2\rho (\delta _1  + \delta _2  + \delta _1 \delta _2 ))T\left\| {a(t)} \right\|_{C[0,T]} \left( {\sum\limits_{k = 1}^\infty  {(\lambda _k^3 \left\| {u_{2k - 1} (t)} \right\|_{C[0,T]} )^2 } } \right)^{\frac{1}{2}} } \right]} \right\},
\end{equation}
where
\[
\rho _1  = \rho ^2 \left[ {\delta _1 \left( {\frac{{\delta _2 }}{2} + \frac{T}{2}} \right) + \delta _2 \left( {\frac{1}{2} + \frac{T}{2} + \delta _2 \left( {\frac{1}{4} + \frac{T}{2}} \right)} \right)} \right.
\]
\[
\left. { + \delta _1 \delta _2 \left( {\frac{1}{2} + \frac{T}{2} + \delta _2 \left( {\frac{1}{2} + \frac{T}{2}} \right)} \right)} \right] + \frac{T}{2} + \delta _2 \left( {\frac{1}{2} + \frac{T}{2}} \right),
\]
\[
\rho _2  = \rho ^2 \left[ {\delta _1 \left( {\frac{1}{2} + \frac{T}{2}} \right) + \delta _1 \left( {\frac{1}{4} + \frac{T}{2}} \right) + \delta _2 \left( {\frac{{\delta _1 }}{2} + \frac{T}{2}} \right)} \right.
\]
\[
\left. { + \delta _1 \delta _2 \left( {\frac{1}{2} + \frac{T}{2} + \delta _1 \left( {\frac{1}{2} + \frac{T}{2}} \right)} \right)} \right] + \frac{1}{2} + \frac{T}{2} + \delta _1 \left( {\frac{1}{2} + \frac{T}{2}} \right).
\]

Then from \eqref{E:23}, \eqref{E:24}, \eqref{E:25}, and \eqref{E:26}, taking into account \eqref{E:10}, \eqref{E:11}, respectively, we find:
\[
\left\| {\tilde u_0 (t)} \right\|_{C[0,T]}\le \frac{2}{{1 + \delta _1 }}\left\| {\varphi (x)} \right\|_{L_2 (0,1)}  + \frac{{2T}}{{1 + \delta _2 }}\left\| {\psi (x)} \right\|_{L_2 (0,1)}  + \frac{{1 + 3\delta _1  + 3\delta _2 }}{{(1 + \delta _1 )(1 + \delta _2 )}}T
\]
\[
 \times (2\sqrt T \left\| {f(x,t)} \right\|_{L_2 (D_T )}  + T\left\| {a(t)} \right\|_{C[0,T]} \left\| {u(x,t)} \right\|_{B_{2,T}^3 } ),
\]
\[
\left\{ {\sum\limits_{k = 1}^\infty  {(\lambda _k^3 \left\| {\tilde u_{2k - 1} (t)} \right\|_{C[0,T]} )^2 } } \right\}^{\frac{1}{2}}
\]
\[\le 4\sqrt 2 \rho (1 + \delta _2 )\left\| {\varphi '''(x)} \right\|_{L_2 (0,1)}+ 4\sqrt 2 \rho (1 + \delta _1 )\left\| {\psi ''(x)} \right\|_{L_2 (0,1)}
\]
\[
 + 4(1 + 2\rho (\delta _1  + \delta _2  + \delta _1 \delta _2 ))\sqrt {2T} \left\| {f_{xx} (x,t)} \right\|_{L_2 (D_T )}
\]
\[
 + 2(1 + 2\rho (\delta _1  + \delta _2  + \delta _1 \delta _2 ))T\left\| {a(t)} \right\|_{C[0,T]} \left\| {u(x,t)} \right\|_{B_{2,T}^3 } ,
\]
\[
\left\{ {\sum\limits_{k = 1}^\infty  {(\lambda _k^3 \left\| {\tilde u_{2k} (t)} \right\|_{C[0,T]} )^2 } } \right\}^{\frac{1}{2}}
\]
\[
\le 8\rho (1 + \delta _2 )\left\| {\varphi '''(x)(1 - q - px) - 3p\varphi ''(x)} \right\|_{L_2 (0,1)}
\]
\[
 + 8\rho (1 + \delta _1 )\left\| {\psi ''(x)(1 - q - px) - 2p\psi '(x)} \right\|_{L_2 (0,1)}
\]
\[
 + 8(1 + 2\rho (\delta _1  + \delta _2  + \delta _1 \delta _2 ))\sqrt T \left\| {f_{xx} (x,t)(1 - q - px) - 2pf_x (x,t)} \right\|_{L_2 (D_T )}
\]
\[
+ 2\sqrt 2 (1 + 2\rho (\delta _1  + \delta _2  + \delta _1 \delta _2 ))T\left\| {a(t)} \right\|_{C[0,T]} \left\| {u(x,t)} \right\|_{B_{2,T}^3 }
\]
\[
+ 8\rho _1 \left\| {\varphi '''(x)} \right\|_{L_2 (0,1)}  + 8\rho _2 \left\| {\psi ''(x)} \right\|_{L_2 (0,1)}
\]
\[
+ 2\sqrt 2 (1 + 2\rho (\delta _1  + \delta _2  + \delta _1 \delta _2 ))^2 T[2\sqrt {2T} \left\| {f_{xx} (x,t)} \right\|_{L_2 (D_T )}
\]
\[
+ T\left\| {a(t)} \right\|_{C[0,T]} \left\| {u(x,t)} \right\|_{B_{2,T}^3 } ],
\]
\[
\left\| {\tilde a(t)} \right\|_{C[0,T]}  \le \left\| {[h(t)]^{ - 1} } \right\|_{C[0,T]}
\]
\[
\times\left\{ {\left\| {h''(t) - f\left( {\frac{1}{2},t} \right)} \right\|_{C[0,T]} } \right. + \frac{1}{2}\left( {\sum\limits_{k = 1}^\infty  {\lambda _k^{ - 2} } } \right)^{\frac{1}{2}} [2\sqrt 2 \rho (1 + \delta _2 )\left\| {\varphi '''(x)} \right\|_{L_2 (0,1)}
\]
\[
 + 2\sqrt 2 \rho (1 + \delta _1 )\left\| {\psi ''(x)} \right\|_{L_2 (0,1)}  + (1 + 2\rho (\delta _1  + \delta _2  + \delta _1 \delta _2 ))2\sqrt {2T} \left\| {f_{xx} (x,t)} \right\|_{L_2 (D_T )}
\]
\[
 + (1 + 2\rho (\delta _1  + \delta _2  + \delta _1 \delta _2 ))T\left\| {a(t)} \right\|_{C[0,T]} \left\| {u(x,t)} \right\|_{B_{2,T}^3 } ]\} ,
\]
or
\begin{equation}\label{E:27}
\left\| {\tilde u_0 (t)} \right\|_{C[0,T]}  \le A_1 (T) + B_1 (T)\left\| {a(t)} \right\|_{C[0,T]} \left\| {u(x,t)} \right\|_{B_{2,T}^3 } ,
\end{equation}
\begin{equation}\label{E:28}
\left\{ {\sum\limits_{k = 1}^\infty  {(\lambda _k^3 \left\| {\tilde u_{2k - 1} (t)} \right\|_{C[0,T]} )^2 } } \right\}^{\frac{1}{2}}\le A_2 (T) + B_2 (T)\left\| {a(t)} \right\|_{C[0,T]} \left\| {u(x,t)} \right\|_{B_{2,T}^3 } ,
\end{equation}
\begin{equation}\label{E:29}
\left\{ {\sum\limits_{k = 1}^\infty  {(\lambda _k^3 \left\| {\tilde u_{2k} (t)} \right\|_{C[0,T]} )^2 } } \right\}^{\frac{1}{2}}\le A_3 (T) + B_3 (T)\left\| {a(t)} \right\|_{C[0,T]} \left\| {u(x,t)} \right\|_{B_{2,T}^3 } ,
\end{equation}
\begin{equation}\label{E:30}
\left\| {\tilde a(t)} \right\|_{C[0,T]}  \le A_4 (T) + B_4 (T)\left\| {a(t)} \right\|_{C[0,T]} \left\| {u(x,t)} \right\|_{B_{2,T}^3 } ,
\end{equation}
where
\[
A_1 (T) = \frac{2}{{1 + \delta _1 }}\left\| {\varphi (x)} \right\|_{L_2 (0,1)}
\]
\[
+ \frac{{2T}}{{1 + \delta _2 }}\left\| {\psi (x)} \right\|_{L_2 (0,1)}  + \frac{{2(1 + 3\delta _1  + 3\delta _2 )}}{{(1 + \delta _1 )(1 + \delta _2 )}}T\sqrt T \left\| {f(x,t)} \right\|_{L_2 (D_T )} ,
\]
\[
B_1 (T) = \frac{{1 + 3\delta _1  + 3\delta _2 }}{{(1 + \delta _1 )(1 + \delta _2 )}}T^2 ,
\]
\[
A_2 (T) = 4\sqrt 2 \rho (1 + \delta _2 )\left\| {\varphi '''(x)} \right\|_{L_2 (0,1)}  + 4\sqrt 2 \rho (1 + \delta _1 )\left\| {\psi ''(x)} \right\|_{L_2 (0,1)}
\]
\[
 + 4(1 + 2\rho (\delta _1  + \delta _2  + \delta _1 \delta _2 ))\sqrt {2T} \left\| {f_{xx} (x,t)} \right\|_{L_2 (D_T )} ,
\]
\[
B_2 (T) = 2(1 + 2\rho (\delta _1  + \delta _2  + \delta _1 \delta _2 ))T,
\]
\[
A_3 (T) = 8\rho (1 + \delta _2 )\left\| {\varphi '''(x)(1 - q - px) - 3p\varphi ''(x)} \right\|_{L_2 (0,1)}
\]
\[
+ 8\rho (1 + \delta _1 )\left\| {\psi ''(x)(1 - q - px) - 2p\psi '(x)} \right\|_{L_2 (0,1)}
\]
\[
 + 8(1 + 2\rho (\delta _1  + \delta _2  + \delta _1 \delta _2 ))\sqrt T \left\| {f_{xx} (x,t)(1 - q - px) - 2pf_x (x,t)} \right\|_{L_2 (D_T )}
\]
\[
+ 8\rho _1 \left\| {\varphi '''(x)} \right\|_{L_2 (0,1)} + 8\rho _2 \left\| {\psi ''(x)} \right\|_{L_2 (0,1)} 
\]
\[
+ 8(1 + 2\rho (\delta _1  + \delta _2  + \delta _1 \delta _2 ))^2 T\sqrt T \left\| {f_{xx} (x,t)} \right\|_{L_2 (D_T )} ,
\]
\[
B_3 (T) = 2\sqrt 2 (1 + 2\rho (\delta _1  + \delta _2  + \delta _1 \delta _2 ))T + 2\sqrt 2 (1 + 2\rho (\delta _1  + \delta _2  + \delta _1 \delta _2 ))^2 T^2 ,
\]
\[
A_4 (T) = \left\| {[h(t)]^{ - 1} } \right\|_{C[0,T]}
\]
\[
 \times \left\{ {\left\| {h''(t) - f\left( {\frac{1}{2},t} \right)} \right\|_{C[0,T]}  + } \right.\frac{1}{2}\left( {\sum\limits_{k = 1}^\infty  {\lambda _k^{ - 2} } } \right)^{\frac{1}{2}} [2\sqrt 2 \rho (1 + \delta _2 )\left\| {\varphi '''(x)} \right\|_{L_2 (0,1)} \]
\[
 + 2\sqrt 2 \rho (1 + \delta _1 )\left\| {\psi ''(x)} \right\|_{L_2 (0,1)}  + (1 + 2\rho (\delta _1  + \delta _2  + \delta _1 \delta _2 ))2\sqrt {2T} \left\| {f_{xx} (x,t)} \right\|_{L_2 (D_T )} ]\} ,
\]
\[
B_4 (T) = \frac{1}{2}\left\| {[h(t)]^{ - 1} } \right\|_{C[0,T]} \left( {\sum\limits_{k = 1}^\infty  {\lambda _k^{ - 2} } } \right)^{\frac{1}{2}} (1 + 2\rho (\delta _1  + \delta _2  + \delta _1 \delta _2 ))T.
\]

Now, from \eqref{E:27}--\eqref{E:29} we obtain:
\begin{equation}\label{E:31}
\left\| {\tilde u(x,t)} \right\|_{B_{2,T}^3 }  \le A_5 (T) + B_5 (T)\left\| {a(t)} \right\|_{C[0,T]} \left\| {u(x,t)} \right\|_{B_{2,T}^3 } ,
\end{equation}
where
\[
A_5 (T) = A_1 (T) + A_2 (T) + A_3 (T), \ \ B_5 (T) = B_1 (T) + B_2 (T) + B_3 (T).
\]

Finally, from \eqref{E:30} and \eqref{E:31} we conclude:
\begin{equation}\label{E:32}
\left\| {\tilde u(x,t)} \right\|_{B_{2,T}^3 }  + \left\| {\tilde a(t)} \right\|_{C[0,T]}  \le A(T) + B(T)\left\| {a(t)} \right\|_{C[0,T]} \left\| {u(x,t)} \right\|_{B_{2,T}^3 } ,
\end{equation}
where
\[
A(T) = A_4 (T) + A_5 (T), \ \ B(T) = B_4 (T) + B_5 (T).
\]

So, let us prove the following theorem.

\begin{theorem}\label{thm:2}
Let conditions $C_1$ - $C_5$ be satisfied and the inequality
\begin{equation}\label{E:33}
B(T)(A(T) + 2)^2  < 1
\end{equation}
holds. Then problem \eqref{E:1}--\eqref{E:3}, \eqref{E:6}, \eqref{E:7} has a unique solution in the ball $K = K_R (\left\| z \right\|_{E_T^3 }  \le R \le A(T) + 2)$ of the space $E_T^3$.
\end{theorem}

\begin{remark}
Inequality \eqref{E:33} is satisfied for sufficiently small values of $T$.
\end{remark}

\begin{proof}
Let us denote $z = [u(x,t),a(t)]^T$ and rewrite the system of equations \eqref{E:20} and \eqref{E:22} in the following operator equation
\begin{equation}\label{E:34}
z = \Phi z,
\end{equation}
where $\Phi  = [\varphi _1 ,\varphi _2 ]^T$ and $\varphi _1 (z),\varphi _2 (z)$ defined by the right sides of \eqref{E:20} and \eqref{E:22}, respectively.

Consider the operator $\Phi (u,a)$ in the ball $K = K_R$ of the space $E_T^3$. Similar to \eqref{E:32} we obtain that for any $z_1 ,z_2 ,z \in K_R$ the following estimates hold:
\begin{equation}\label{E:35}
\left\| {\Phi z} \right\|_{E_T^3 }  \le A(T) + B(T)\left\| {a(t)} \right\|_{C[0,T]} \left\| {u(x,t)} \right\|_{B_{2,T}^3 }\le A(T) + B(T)(A(T) + 2)^2 ,
\end{equation}
\begin{equation}\label{E:36}
\left\| {\Phi z_1  - \Phi z_2 } \right\|_{E_T^3 } \le B(T)R(\left\| {a_1 (t) - a_2 (t)} \right\|_{C[0,T]}  + \left\| {u_1 (x,t) - u_2 (x,t)} \right\|_{B_{2,T}^3 } ).
\end{equation}

Then taking into account \eqref{E:33} in \eqref{E:35} and \eqref{E:36}, it follows that the operator $\Phi$ acts in the ball $K = K_R$ and is contractive. Therefore, in the ball $K = K_R$ the operator $\Phi$ has a unique fixed point $\{ z\}  = \{ u,a\}$ that is a unique solution of \eqref{E:34} in the ball $K = K_R$; i.e. it is a unique solution of system \eqref{E:20}, \eqref{E:22} in the ball $K = K_R$.
\end{proof}

Thus, we obtain that the function $u(x,t)$ as an element of the space $B_{2,T}^3$ is continuous and has continuous derivatives $u_x (x,t)$ and $u_{xx} (x,t)$ in $D_T$.

Analogous to \cite{MehYus} one can show that the function $u_{tt} (x,t)$ is continuous in the region $D_T$.

It is easy to verify that Eq.\eqref{E:1} and conditions \eqref{E:2}, \eqref{E:3}, \eqref{E:6}, and \eqref{E:7} are satisfied in the ordinary sense. Consequently, $\{ u(x,t),a(t)\}$ is a solution of problem \eqref{E:1}--\eqref{E:3}, \eqref{E:6}, \eqref{E:7} and by Remark 3.1 this solution is unique in the ball $K = K_R$. The theorem is proved.

From Theorem 1 and Theorem 2 immediately imply that the original problem \eqref{E:1}--\eqref{E:5} has a unique classical solution.

\begin{theorem}\label{thm:3}
Suppose that all the conditions of Theorem 2 are satisfied and
\[
\int\limits_0^1 {f(x,t)dx}  = 0, \ 0 \le t \le T, \ \  \int\limits_0^1 {\varphi (x)dx}  = 0, \ \ \int\limits_0^1 {\psi (x)dx}  = 0,
\]
\[
h(0) + \delta _1 h(T) = \varphi \left( {\frac{1}{2}} \right), \ \ h'(0) + \delta _2 h'(T) = \psi \left( {\frac{1}{2}} \right),
\]
\[
\frac{{(1 + 2\delta _1  + 3\delta _2  + \delta _1 \delta _2 )T^2 (A(T) + 2)}}{{2(1 + \delta _1 )(1 + \delta _2 )}} < 1.
\]

Then problem \eqref{E:1}--\eqref{E:5} has a unique classical solution in the ball $K = K_R (\left\| z \right\|_{E_T^3 } \\ \le R \le A(T) + 2)$ of the space $E_T^3$.
\end{theorem}

\section{Conclusion}\label{Sec:Conclusion}
In the work, the classical solvability of a nonlinear inverse boundary value problem for a second-order hyperbolic equation with nonlocal conditions was studied. The considered problem was reduced to an auxiliary inverse boundary value problem in a certain sense and its equivalence to the original problem is shown. Then using the Fourier method and contraction mappings principle, the existence and uniqueness theorem for auxiliary problem is proved. Further, on the basis of the equivalency of these problems, the existence and uniqueness theorem for the classical solution of the original inverse coefficient problem is established for the smaller value of time.

% A simple way to write an acknowledgement
\section*{Acknowledgements}
The authors express their profound thanks to the anonymous reviewers for their insightful comments and helpful suggestions.

\end{document}